\documentclass{article}
\usepackage{amsmath,amsthm,verbatim,amssymb,amsfonts,amscd,diagrams, graphics}
\theoremstyle{plain}
\newtheorem{theorem}{Theorem}
\newtheorem{lemma}[theorem]{Lemma}
\newtheorem{proposition}[theorem]{Proposition}

\theoremstyle{definition}

\theoremstyle{remark}
\newtheorem{remark}[theorem]{Remark}

\newarrow{ul}---->
\newarrow{Backwards}<----

\newcommand{\td}[1]{\tilde{#1}}
\newcommand{\into}{\hookrightarrow}
\newcommand{\Z}{\mathbb{Z}}

\newcommand{\R}{\mathbb{R}}

\newcommand{\bd}{\partial}

\begin{document}

\title{Some non-trivial PL knots whose complements\\are homotopy circles}
\author{Greg Friedman\\ Yale University}
\date{July 7, 2005 }

\maketitle
%typeset=\today
%\tableofcontents



\begin{abstract}
We show that there exist non-trivial piecewise-linear (PL) knots with isolated singularities $S^{n-2}\subset  S^n$, $n\geq 5$, whose complements have the homotopy type of a circle. This is in contrast to the case of smooth, PL locally-flat, and topological locally-flat  knots,  for which it is known that if the complement has the homotopy type of a circle, then the knot is trivial. 
\end{abstract}

It is well-known that if the complement of a smooth, piecewise linear (PL) locally-flat, or topological locally-flat  knot $K \subset S^n$, $K\cong S^{n-2}$, $n\geq 5$, has the homotopy type of a circle, then $K$ is equivalent to the standard unknot in the appropriate category (see Stallings \cite{St63} for the topological case and Levine \cite{L65} and \cite[\S 23]{L70} for the smooth and PL cases). This is also true of classical knots $S^1\into S^3$ (see \cite[\S 4.B]{Rolf}), for which these categories are all equivalent, and in the topological category for knots $S^2\into S^4$ by Freedman \cite[Theorem 6]{Fr84}. 

By contrast, Freedman and Quinn showed in \cite[\S 11.7]{FrQ} that any classical knot with Alexander polynomial $1$ bounds a topological locally-flat $D^2$ in $D^4$ whose complement is a homotopy circle, and by collapsing the boundary, one obtains a singular $S^2$ in $S^4$ with the same property. In the same dimensions, Boersema and Taylor \cite{BT} constructed a specific example of a PL knot with an isolated singularity whose complement is a homotopy circle. It follows by taking iterated suspensions that there are PL knots in all dimensions $n\geq 4$ whose complements are homotopy circles, though this process will lead to increasingly more complicated singularities. 
In this note, we construct PL knots for any $n\geq 5$ that are locally-flat except at one point and whose complements are homotopy circles. 

To construct the  knots with the desired properties,
it will suffice to construct for each $n\geq 5$ a  PL locally-flat  disk knot $L \subset D^n$, such that $D^n-L\sim_{h.e.} S^1$ and such that the PL locally-flat boundary sphere knot $\bd L\subset \bd D^n$ is non-trivial. 
By a PL locally-flat disk knot $L\subset D^n$, we  mean the image of a PL locally-flat embedding $D^{n-2}\into D^n$ such that $\bd L\subset \bd D^n$ is a locally-flat sphere knot and $\text{\emph{int}} (L)\subset \text{\emph{int}}(D^n)$. This will suffice since, if such a disk knot exists, we may then adjoin the cone on the boundary pair $(\bd D^n, \bd L)$ to obtain  a PL sphere knot $K\subset S^n$ that is locally-flat except at the cone point: 
\begin{diagram}
K&=& L& \cup_{\bd L}& c(\bd L) \\
\cap& &\cap&&\cap\\
S^n&=& D^n&\cup_{\bd D^n} &c(\bd D^n)&.
\end{diagram}
It is clear that $S^n-K\sim_{h.e.} D^n-L$, so  if the complement of $L$ is a homotopy circle then so will be that of $K$. Furthermore, $K$ will be non-trivial since the link pair of the cone point will be non-trivially knotted, which is impossible in the unknot, which is locally-flat.

So we construct such a disk knot. The procedure will be based upon that given by the author in \cite{GBF1} for constructing certain Alexander polynomials of disk knots, which in turn was a generalization of Levine's construction of sphere knots with given Alexander polynomials in \cite{L66}. All spaces and maps will be in the PL category without further explicit mention. 

Suppose that $n\geq 5$, and let $U$ be the trivial disk knot $U\subset  D^n$, i.e. $D^n$ may be identified with the unit ball in $\R^n$ such that $U$ is the intersection of $D^n$ with the coordinate plane $\R^{n-2}\subset \R^n$. Embed an unknotted  $S^{n-3}$ into $\bd D^n=S^{n-1}$  so that it is not linked with $\bd U$ (in fact, we may assume that the new $S^{n-3}$ and $\bd U$ are in opposite  hemispheres of $\bd D^n$). We use the standard framing of the new unknotted $S^{n-3}$ to attach an $n-2$ handle to $D^n$, obtaining a space homeomorphic to $S^{n-2}\times D^2$ and containing an unknotted disk in a trivial neighborhood of some point on the boundary. We can assume that $U$ bounds an  embedded $n-1$ disk $V$ in $S^{n-2}\times D^2$, that $\bd U$ bounds an $n-2$ disk $F$ in $\bd (S^{n-2}\times D^2)$, that $\bd V=U\cup F$, and that  $\text{\emph{int}}(V)\subset \text{\emph{int}}(S^{n-2}\times D^2)$.  Let $C_0=S^{n-2}\times D^2-U$, and let $\td C_0$ be the infinite cyclic cover of $C_0$ associated with the kernel of the homomorphism $\pi_1(C_0)=\Z\to \Z$ determined by linking number with $U$. Let $X_0=\bd (S^{n-2}\times D^2)-\bd U$, and let $\td X_0$ be the infinite cyclic cover of $X_0$ in $\td C_0$.

As in the usual construction of infinite cyclic covers in knot theory (see, e.g., Rolfsen \cite{Rolf}), we can form $\td C_0$ by a cut and paste procedure: we cut $C_0$ along $V$ to obtain $Y_0$ and then glue a countably infinite number of copies of $Y_0$ together along the copies of $V$. Since $C_0-V\sim_{h.e.}S^{n-2}$, we have $\td H_{n-2}(\td C_0)=\Z[\Z]=\Z[t,t^{-1}]$ - where $t$ represents a generator of the group of covering translations - and all other reduced homology groups are trivial. Similarly, since  $\bd (S^{n-2}\times D^2)-F$ is a punctured $S^{n-2}\times S^1$, $\td H_*(\td X_0)$ is $\Z[\Z]$ in dimensions $n-2$ and $1$, and trivial otherwise. 

It is also apparent  that $\pi_*(\td C_0)$ is trivial for $*<n-2$, while $\pi_1(\td X_{0})$ is free on a countably infinite number of generators. Thus,  since $n\geq 5$, $\pi_2(\td C_0,\td X_0)$ is also  free on a countably infinite number of generators. Meanwhile, for $X_0$, itself, $\pi_1(X_0)$ is the free group on two generators: one generator corresponds to the generator of $\pi_1(\bd (S^{n-2}\times D^2))=\pi_1(S^{n-2}\times S^1)=\Z$ and the other corresponds to the meridian of the unknotted $\bd U$ (this can be demonstrated by an easy Seifert-van Kampen argument, by considering $\bd U$ to lie in a ball neighborhood of some point). Let $a$ represent the generator corresponding to the meridian of $\bd U$, and let $b$ represent the other described generator. Similarly, $\pi_1(C_0)\cong \Z$, its generator also being given by $a$, while $b$ is contractible  in this larger space. 

Consider now the element $\gamma$ of $\pi_1(X_0)$ given by $b^2aba^{-1}b^{-1}ab^{-1}a^{-1}$. Since $b=1$  in $\pi_1(C_0)$ and $a$ occurs with total exponent $0$ in $\gamma$, the image of $\gamma$ in $\pi_1(C_0)$ is trivial, so  any representative of $\gamma$ is the boundary of a $2$-disk $\Gamma$ in $C_0$. Since $n\geq 5$, we can assume that $\Gamma$ is properly embedded (see \cite[Corollary 8.2.1]{Hud}). Furthermore, $\gamma$ can be lifted to a closed curve in $\td X_0$; if we let $c_i$ represent the generators of $\pi_1(\td X_0)$, then any lift of $a$ is a path between adjoining lifts of $X_0$ in the cut and paste construction, and $\gamma$ lifts to $\td \gamma=c_0^2c_1c_0^{-1}c_1^{-1}\in \pi_1(\td X_0)$. In the abelianization $H_1(\td X_0)$, the image of $\td \gamma$ is the same as the image of $c_0$, which is a $\Z[\Z]$-module generator of $H_1(\td X_0)$.

Let $N$ denote an open regular neighborhood of $\Gamma$ in $C_0$.  We claim that $S^{n-2}\times D^2-N$ is homeomorphic to $D^n$. In fact, observe that in $S^{n-2}\times S^1$, $\gamma$ is homotopic to the standard generator $b=*\times S^1$ of $\pi_1( S^{n-2}\times S^1 )$ (with an appropriate choice of orientations). Thus, in $(S^{n-2}\times D^2, S^{n-2}\times S^1)$, the pair $(\Gamma, \gamma)$ is homotopic to the standard generator $*\times D^2$ of $\pi_2(S^{n-2}\times D^2,S^{n-2}\times S^1)$. These homotopies can be realized by ambient isotopies by  \cite[Theorem 10.2]{Hud}.  Then it is clear that $ S^{n-2}\times D^2 -N\cong D^{n-2}\times D^2\cong D^n$.

\begin{comment}Then it is clear that $ S^{n-2}\times D^2 -N\cong D^{n-2}\times D^2\cong D^n$. If $n=5$, the result still holds, but the the argument is slightly more complicated. Full details can be found in \cite[pp. 1512-3]{GBF1}, but the idea is as follows: By dimensional considerations, there is still an ambient isotopy of $S^{n-2}\times S^1$ that takes $\gamma$ to $*\times S^1$; from this it can be seen that $\bd S^{n-2}\times D^2-N\cong S^{n-1}$. It then suffices to show that $\Delta=\bd S^{n-2}\times D^2-N$ is contractible, since it follows from the Poincar\'e conjecture that any contractible $5$-manifold whose boundary is a sphere is homeomorphic to $D^5$. Direct calculations can then be used to show that $\pi_1(\Delta)=1$ and $\td H_*(\Delta)=0$, from which the claim follows by the Hurewicz Theorem. 
\end{comment}

Fixing a homeomorphism $S^{n-2}\times D^2 -N\to  D^n$,  the  image of $U$ is a new disk knot, which we christen $L$. We claim that $L$ is no longer trivial but that its complement is a homotopy circle. 

Let $C$ be the complement of an open regular neighborhood of $L$ in $D^n$ (the disk knot exterior). Thus $C$ is homotopy equivalent to $D^n-L$. Similarly, let $X$ be the exterior of $\bd L$ in $\bd D^n=S^{n-1}$. We must study the homotopy and homology of $C$, $X$, and their coverings.

\begin{lemma}\label{L: 1} $\pi_1(C)=\Z$.
\end{lemma}
\begin{proof} $C\sim_{h.e.} D^n-L\cong C_0-N$, and since $N$ is the regular neighborhood of $\Gamma\cong D^2$, $N$ is homeomorphic to $D^n$ and $\bd N\cong D^2\times S^{n-3}$. So, up to homeomorphism, we may think of $C_0$ as $(C^0-N)\cup_{D^n\times S^{n-3}} D^n$. Since $n\geq 5$, we see from the Seifert-van Kampen Theorem that $\pi_1(C_0-N)\cong \pi_1(C_0)$. Since $\pi_1(\td C_0)=1$, where $\td C_0$ is the infinite cyclic cover of $C_0$, it follows that $\pi_1(C_0)\cong \Z$. Thus $\pi_1(C)\cong \Z$. \end{proof}

\begin{lemma} $\pi_1(X)\cong \langle a,b\mid b^2aba^{-1}b^{-1}ab^{-1}a^{-1}\rangle$.\end{lemma}
\begin{proof} The effect of the handle subtraction $C_0-N$ on the boundary $X_0$ is that of a surgery on the embedded curve $\gamma$. Since $\pi_1(X_0)$ is free on the generators $a$ and $b$, the result of the surgery is the given group. 
(Proof: The result of the surgery is $(X_0-S^1\times D^{n-2}) \cup D^2\times S^{n-3}$, where the $S^1$ represents $\gamma$. But since $n\geq 5$, $ \pi_1(X_0-S^1\times D^{n-2})\cong \pi_1(X_0)$. So by Seifert-van Kampen, $\pi_1$ of the result of the surgery is $\pi_1(X_0)/\pi_1(S^1\times S^{n-3})\cong \pi_1(X_0)/\Z$, where the $\Z$ is generated by $S^1\times *$ in $S^1\times S^{n-3}$, which is the boundary of the neighborhood of $\gamma$. But any such curve  is homotopic to $\gamma$, which represents $b^2aba^{-1}b^{-1}ab^{-1}a^{-1}$.)
\end{proof}

\begin{lemma}\label{L: 2} The Alexander modules $\td H_*(\td C)$, $\td H_*(\td X)$, and $\td H_*(\td C,\td X)$  are all trivial.\end{lemma}\begin{proof} Let $\td \gamma$ be the lift of $\gamma$ considered above. We can also lift $\Gamma$ to a $2$-disk $\td \Gamma$ in $\td C_0$. In fact, we can find a countable number of lifts $\td \gamma_i$ and $\td \Gamma_i$, and, since $\Gamma$ is embedded,  the $\td \Gamma_i$ are all disjoint. If $\td N_i$ then represent the lifts of the regular neighborhood $N$, $\td C_0-\amalg_i \td N_i$ will be the infinite cyclic cover of $ C_0-N\cong D^n-L$.  

Now consider $\td X_0\cup \amalg_i\td N_i$. Each intersection $\td X_0\cap N_i$ is homotopy equivalent to a translate of $\td \gamma_i$, which we know represents the $\Z[\Z]$-module generator of $H_1(\td X^0)$. It thus follows from the Mayer-Vietoris sequence that $\td H_*(\td X_0\cup \amalg_i\td N_i)$ is trivial except in dimension $n-2$, where it is $\Z[\Z]$. Meanwhile, we already know that $\td H_*(\td C_0)$ is trivial except in dimension $n-2$, where it is also $\Z[\Z]$. Consider the map $H_{n-2}(\td X_0\cup \amalg_i\td N_i)\to H_*(\td C_0)$. In each module, a $\Z[\Z]$-module generator is represented by a choice of $S^{n-2}\times *\subset S^{n-2}\times S^1\subset S^{n-2}\times D^2$ that is disjoint from $V$. Thus this homology map is an isomorphism, and it follows that $H_*(\td C_0, \td X_0\cup \amalg_i\td N_i)$ is trivial. 
But by excision, $H_*(\td C_0, \td X_0\cup \amalg_i\td N_i)\cong H_*(\td C,\td X)$. 

Similarly, it follows from easy homological calculations that $\td H_*(\td X)$ is trivial. In fact, it can be seen that the construction of $X$ from $X_0$ is by a surgery, and upon restriction of our construction to its effect on $X_0$, we obtain the construction of  Levine for producing smooth sphere knots with given Alexander polynomials in \cite{L66}. In this case, the Alexander polynomial is trivial (since $\td \gamma$ generates $H_1(\td X_0)$), and it follows from Levine's calculations that $\td H_*(\td X)=0$.

Then $\td H_*(\td C)$ is also trivial, by the long exact sequence of the pair $(\td C,\td X)$. 
\end{proof}

\begin{proposition}
$\pi_*(D^n-L)\cong \pi_*(S^1)$.
\end{proposition}
\begin{proof}
By Lemma \ref{L: 1},  $\pi_1(C)=\Z$. Thus the infinite cyclic cover $\td C$ is simply connected, and  
since we also have $\td H_*(\td C)=0$ by Lemma \ref{L: 2}, it follows that $\pi_j(\td C)=0$ for all $j>1$ by Hurewicz's Theorem. Thus for $j>1$, $\pi_j(C)=0$, and $\pi_*(D^n-L)\cong \pi_*(C)\cong \pi_*(S^1)$.
\end{proof}

\begin{theorem}
$D^n-L$ is a homotopy circle.
\end{theorem}
\begin{proof}
By the preceding proposition, $D^n-L$ has the same homotopy groups as a circle. But $D^n-L$ is homotopy equivalent to $C$, which is homeomorphic to a finite simplicial complex. Since the inclusion $i:S^1\to C$ of a meridian of $L$ induces the isomorphism $\pi_1(S^1)\to \pi_1(C)$, we can conclude that $i$ is a homotopy equivalence. Thus 
 $C\sim_{h.e.} D^n-L$ is a homotopy circle. 
\end{proof}

It only remains to show that $L$ is non-trivial, which will follow once we show that the group $\pi_1(X)$ of the boundary knot $\bd L$ is not $\Z$. 

\begin{lemma}
The group $G=\langle a,b\mid b^2aba^{-1}b^{-1}ab^{-1}a^{-1}\rangle$ is not isomorphic to $\Z$.
\end{lemma}
\begin{proof}
This lemma can be proven in a variety of ways. The following elegant demonstration was shown to me by Andrew Casson.

We adjoin an extra generator $c$, which we immediately set equal to $aba^{-1}$. Then
\begin{align*}
\langle a,b\mid b^2aba^{-1}b^{-1}ab^{-1}a^{-1}\rangle&\cong \langle a,b,c \mid b^2aba^{-1}b^{-1}ab^{-1}a^{-1}, cab^{-1}a^{-1}\rangle\\
&\cong \langle a,b,c \mid  b^2cb^{-1}c^{-1}, cab^{-1}a^{-1}\rangle\\
&\cong \frac{ \langle b,c \mid  b^2cb^{-1}c^{-1}\rangle * \langle a \rangle}{\langle cab^{-1}a^{-1}\rangle }.
\end{align*}
Written this way, $G$ has the form of an HNN extension of the Baumslag-Solitar group $H=\langle b,c \mid  b^2cb^{-1}c^{-1}\rangle$, which is isomorphic to the semi-direct product $\Z[\frac{1}{2}]\rtimes \Z$. Thus $H$ is a non-abelian subgroup of $G$, which hence cannot be $\Z$. 

Alternatively, to apply an unnecessarily large hammer, once $G$ is written as $\langle a,b,c \mid b^2aba^{-1}b^{-1}ab^{-1}a^{-1}, cab^{-1}a^{-1}\rangle$, it follows from \cite{RTV} that $G$ is not even residually finite. 

A third proof would utilize Whitehead's theorem on one-relator groups \cite{Wh36}. 
\end{proof}

\begin{remark}
There is nothing exceptionally special about the group $G$ we have used in this construction, except that it turned out to be a fairly tractable example of a group with suitable properties. Any group possessing a two generator, one relator presentation with the properties employed above clearly would be sufficient.  
\end{remark}

\begin{comment}
\begin{itemize}
\item $G\ncong \Z$,
\item $G$ abelianizes to $\Z$,
\item the commutator subgroup of $G$ abelianizes to $0$,
\item and $G$ has a presentation in two generators and one relator such that the group is the normal closure of one of the generators $a$ and the total exponent of $a$ in the relator is $0$.  
\end{itemize}
\end{remark}
\end{comment}

\bibliographystyle{amsplain}
\bibliography{bib}

\end{document}